\documentclass{birkjour_t2}
\usepackage{lineno,hyperref}
\usepackage[utf8]{inputenc}
\usepackage{xcolor}
\usepackage{graphicx}%
\usepackage{multirow}%
\usepackage{amsmath,amssymb,amsfonts}%
\usepackage{amsthm}%
\usepackage{mathrsfs}%
\usepackage[title]{appendix}%
\usepackage{xcolor}%
\usepackage{textcomp}%
\usepackage{manyfoot}%
\usepackage{booktabs}%
\usepackage{algorithm}%
\usepackage{algorithmicx}%
\usepackage{algpseudocode}%
\usepackage{listings}%
\newcommand{\beqa}{\begin{eqnarray}}
\newcommand{\beqan}{\begin{eqnarray*}}
\newcommand{\eeqa}{\end{eqnarray}}
\newcommand{\eeqan}{\end{eqnarray*}}

\theoremstyle{thmstyleone}%
\newtheorem{theorem}{Theorem}
\newtheorem{Lem}{Lemma}
\theoremstyle{thmstyletwo}%
\usepackage{color}
\newcommand{\R}{\mathbb{R}}
\begin{document}

%
%
%
%
%
%
%
%
%

\title[Exponential Stability of Boundary-damped Wave Equation]
{\begin{center}
Revisiting the Direct Fourier Filtering Technique for the Maximal Decay Rate of Boundary-damped Wave Equation by Finite Differences and Finite Elements
\end{center}}



\author[Ozer]{A. \"{O}. \"{O}zer}
\address{Department of Mathematics, Western Kentucky University, \\ Bowling Green, KY 42101, USA. \\
E-mail: {\tt ozkan.ozer@wku.edu}}


\author[Ozer]{Rafi Emran}
\address{Department of Mathematics, Western Kentucky University, \\ Bowling Green, KY 42101, USA. \\
E-mail: {\tt mdrafias.sadeqibnemran227@topper.wku.edu}}




\date{January 1, 2004}

\begin{abstract}
The  one-dimensional PDE model of the  wave equation with a  state feedback controller at its boundary, which describes wave dynamics of a wide-range of controlled mechanical systems,   has exponentially stable solutions. However, it is known that the reduced models of the wave equation by the  standard Finite Differences and Finite Elements suffer from the lack of exponential stability (and exact observability without a state feedback controller) uniformly as the discretization parameter tends to zero.  This is due to the loss of uniform gap  among the high-frequency eigenvalues as the discretization parameter tends to zero. One common remedy to overcome this discrepancy is the direct Fourier filtering of the reduced models, where the high-frequency spurious eigenvalues are filtered out. After filtering, besides from the strong convergency, the exponential decay rate, mimicking the one for the partial differential equation counterpart, can be retained uniformly. However, the existing results in the literature are solely based on an observability inequality of the control-free model, to which the filtering is implemented. Moreover, the decay rate as a function of the filtering parameter is implicit. In this paper, exponential stability results for both filtered Finite Difference and Finite Element reduced models are established directly by a Lyapunov-based approach and a thorough eigenvalue estimation.The maximal decay rate  is explicitly provided as a function of the feedback gain and filtering parameter.  Our results, expectedly, mimic the ones of the PDE counterpart uniformly as the discretization parameter tends to zero. Several numerical tests are provided to support our results.
\end{abstract}

\maketitle

\section{Introduction}\label{sec1}

Let dots denote the derivatives with respect to the time variable $t$. Consider the standard one-dimensional  wave equation clamped at the left. The boundary feedback is injected at the free end:
\begin{eqnarray}
\label{estatic}&&\left\{
  \begin{array}{ll}
 \ddot v-c^2 v_{xx}=0, \quad \left(x,t\right)\in \left(0,L\right)\times\mathbb{R}^+&\\
v\left(0,t\right)=0,  \quad c^2 v_x  \left(L,t\right)=-\xi v_t\left(L,t\right),~~ t\in\mathbb{R}^+&\\
   \left[v,v_t\right]\left(x,0\right)=\left[v_0,v_1\right]\left(x\right),  \quad x\in\left[0,L\right]
   \end{array}\right.
\end{eqnarray}
where $c>0$ is the wave propagation speed and $\xi>0$ is the feedback gain. The closed-loop system \eqref{estatic} is widely used in the literature  to demonstrate transverse vibrations on a  string \cite{L-L}, longitudinal vibrations on a piezoelectric beam \cite{M-O}, or sound vibrations in a duct \cite{Morris}, etc. 

The natural energy of the solutions for \eqref{estatic}  is defined as
\begin{eqnarray}\label{eq4}
       E(t)=\frac{1}{2}\int^L_0\left[\left|\dot v\right|^2+c^2 |v_x|^2\right] dx.
\end{eqnarray}
Defining the state $\vec\psi=(\psi^1,\psi^2)^{\rm T}=(v,\dot v)^{\rm T},$  the system \eqref{estatic} can be formed into the first-order form
\begin{equation}\label{PDE}
	\left\{	\begin{array}{ll}
		{\dot{\vec\psi}}=\mathcal A\vec{\psi}=\begin{bmatrix}
  0     & I   \\
 c^2 \frac{\partial^2}{\partial x^2}     & 0  \\
\end{bmatrix} \vec\psi,\quad  t\in \mathbb{R}^+&\\
	\psi^1(0)=0,\quad \psi^1_x(L)=-\xi \psi^2(L),\\
	\vec \psi(x,0)=\vec \psi_0(x).
	\end{array}\right.
\end{equation}

In fact, the eigenvalues $\{\lambda(\xi)\}_{k\in \mathbb{Z}}$ of the non-self-adjoint system operator $\mathcal A$  in \eqref{PDE} are calculated explicitly \cite{Ito} as the following
				\begin{eqnarray}\label{eigPDE}
				\lambda_k(\xi)=
				\begin{array}{ll}
				-\frac{c}{2L} ln\left|\frac{\xi +c}{\xi -c}\right| + \frac{(2k+1)\pi c }{2L}i,&~~  \xi<c.
				\end{array}
				\end{eqnarray}
Therefore,  the fed-back sensor measurement,$\dot v(L,t)$  is sufficient to push all the eigenvalues to the left-half plane, uniformly bounded away from the imaginary axis. This leads to  the known exponential stability result \cite{Ito,Kom}.

The same exponential stability result can be also proved directly by a Lyapunov approach. A Lyapunov function $E_\delta (t)$  is constructed by a $\delta-$perturbation of the energy $E(t)$ by a standard energy-like  functional  $F(t)$  as the following
\begin{eqnarray}\label{dm1}
F(t):=\int_0^L 2 \dot u x u_x dx,\quad E_\delta(t):=E(t)+\delta F(t).
\end{eqnarray}

\begin{theorem}\label{thm1} \cite[Chap. 5]{Book}
        Let  $H^1_L(0,L):=\{v\in H^1(0,L)~:~ v(0)=0\}$.  There exist constants $M,\delta, \omega>0$ such that  for all initial data $(v^0,v^1)\in H^1_L(0,L)\times L^2(0,L)$, the solutions $(v,\dot v)$ of the controlled system \eqref{estatic} with $0<\xi<c$ are exponentially stable, i.e. 
				$$E(v,\dot v;t)\leq M(\delta) E(v_0,v_1;0)e^{-\sigma (\delta(\xi)) t}$$
where  
\begin{eqnarray}\label{Lius}
\begin{array}{ll}
&\delta(\xi)=\frac{1}{2} \min\left\{\frac{c}{2L},\quad \frac{\xi c^2}{L(c^2+\xi^2)}\right\},\quad \sigma(\delta)=2\delta\left(1-\frac{2L\delta}{c}\right), \quad M(\delta)=\frac{c+2L\delta}{c-2L\delta}.
\end{array}
\end{eqnarray}
Hence, the maximal decay rate $\sigma_{max}=\frac{c}{4L}$ is attained with the optimal feedback gain $\xi=c.$ Indeed, at $\xi=c,$ the spectral abscissa approaches $-\infty.$ Therefore, the solutions disappear in finite time.
\end{theorem}

It is also known that the system \eqref{estatic} with $\xi=0$ is exactly observable in the energy space $H^1_L(0,L)\times L^2(0,L)$. 
Indeed, proving Theorem \ref{thm1} is mathematically equivalent to proving the so-called an observability inequality for the control-free problem, i.e. $\xi=0.$
\begin{theorem}\cite[Chap. 3]{Kom}\label{observa}
 Consider the solutions $(v,\dot v)\in H^1_L(0,L)\times L^2(0,L)$ of the control-free system \eqref{estatic}, i.e. $\xi=0$. For  any $T>\frac{2L}{c}$,  there exists a constant $C(T)>0$ such that  \begin{eqnarray}\label{obsvr}
 \int^T_0 \left|\dot v(L,t)\right|^2~dt\geq C(T)E(0)
 \end{eqnarray}
for all initial data  $(v_0,v_1)\in H^1_L(0,L)\times L^2(0,L).$
\end{theorem}

It is widely known that Theorem \ref{thm1} above fails to hold true as the discretization parameter $h$ in Finite Difference and Finite Element-based approximations tends to zero. This discrepancy is first observed  in the pioneer work of \cite{Ito}, where the widely used model reductions failed to mimic the PDE counterpart \eqref{estatic}. The lack of exponential stability of Finite Difference and Finite Element-based model reductions as the discretization parameter approaches zero is shown rigorously in \cite{Peichl} by thorough spectral estimates. Indeed, it is observed that the high-frequency of eigenvalues tend to the imaginary axis as the discretization parameter approaches zero. Later, it is shown by a multipliers approach that the lack exponential stability   is due to the the lack of exact observability of the control-free model, i.e. $\xi\equiv 0$,   as the discretization parameter  in both Finite Differences and Finite Element approximations tend to zero. In other words,  the observability constant $C(T)$ in \eqref{obsvr} of  Theorem \ref{observa} turns out to blow up to infinity as the discretization parameter approaches zero. Investigating this discrepancy further shows that high-frequency eigenvalues tend to a specific value as the discretization parameter approaches zero. To remedy this issue,  a ``direct Fourier filtering technique'', controlling only the low-frequency  part of the solution in order to eliminate the short wave-length (high-frequency) components of the solutions,  is  proposed for the first time in \cite{I-Z} for fully-clamped boundary conditions. With this approach, the observability results is recovered fully, see also the detailed review paper on this issue \cite{Zuazua}.  Since the numerical approximation is proved to converge to the PDE counterpart, the high-frequency components of the solutions are retained by choosing the discretization small enough.

The first attempt of a Finite-Difference-based numerical scheme for \eqref{estatic}  is reported in \cite{T-Z}. The proof of the exponential stability result is solely based on  an observability result,  and   a numerical viscosity term  artificially added to the system \eqref{estatic}, which is referred to the ``indirect filtering'' technique in the literature. The proof  of exponential stability uses the decomposition of the PDE model \eqref{estatic} into a control-free problem with non-zero initial conditions and the controlled problem with zero initial conditions. The techniques in the proof  prevents the rigorous investigation of  the maximal decay rate  in terms  of the numerical filtering parameter. Later on,  with the direct Fourier filtering technique, as in \cite{I-Z}, the exponential stability for the Finite Difference-based model  is shown to be retained \cite{B} as the discretization parameter approaches zero. The proof outlined in \cite{B} is based off of the similar decomposition technique as in \cite{T-Z} yet the ``direct Fourier filtering'' is implemented to the control-free model, $\xi=0$ in \eqref{estatic}. 
The major drawback of  the exponential stability results in both \cite{B} and \cite{T-Z} is that the proof of the exponential stability result solely relies on an observability result of the control-free model. This together with the decomposition argument make the maximal decay rate analysis, and therefore the analysis of finding the optimal feedback gain to achieve the maximal decay rate, more complicated.  

 It is worthwhile to mention that   an alternate semi-discretized Finite-Difference based model reduction,  based on reducing the order of the system \eqref{estatic} in time and space variables,
is reported \cite{Guo3,Guo4}.  This type of model reduction, similar to the one by  mixed-finite elements \cite{Ito,Banks,Cas},  does not need any numerical filtering.  There are several other remedies proposed in the literature for interested audience worth to read for, e.g. Tichonoff regularization \cite{G},  non-uniform meshes \cite{Erv}, etc. 

In this paper, following the work in \cite{Peichl},  the spectral investigation of the non-self-adjoint system operator $\mathcal{A}$ in \eqref{PDE} is extended for the estimation of the maximum modulus of the eigenvalues of the system matrices for each Finite Difference and Finite Element-based model reductions of \eqref{estatic}.  It is proved that both reduced models lack uniform observability with respect to mesh parameter $h\to 0$.  By the implementation of the  direct Fourier filtering technique to the closed-loop model reductions directly, the exponential stability of the model reductions  with respect to mesh parameter is recovered uniformly. The maximal decay rate for each model reduction is established as a function of the filtering parameter. The exponential stability results mimic the ones in Theorem \ref{thm1}, and the proofs are solely based on a Lyapunov approach, discrete multipliers, and through spectral estimates. Maximal decay rates for both approximations as a function of the filtering parameter are also provided for each approximation technique. 

The overall methodology presented here  not only extends the results in \cite{F,B,T-Z} but also provides  better insights to understand the overall  exponential stability of the Fourier-filtered solutions of \eqref{estatic} by Finite Differences and Finite Elements \cite{Ito,Peichl}.  To the best of our knowledge, the exponential stability of the  Finite Elements with direct Fourier filtering is not reported at all. More importantly, our analysis is applicable to large collections of wave and beam models.

\section{Semi-discretizations of \eqref{estatic} in the $x-$variable}
Let $N\in\mathbb{N}$ be given, and define the mesh size $h:=\frac{1}{N+1}$. Consider a uniform discretization of the interval $[0,L]$: $0=x_0<x_1<...<x_{N-1}<x_N<x_{N+1}=L.$

\subsection{Finite Differences with $\xi\equiv 0$}
Let $v_j=v_j(t)\approx v(x_j,t)$ the approximation of the solution $v(x,t)$ of \eqref{estatic} at the point space $x_j=j\cdot h$ for any $j=0,1,...,N,N+1$, and $\vec{v}=[v_1,v_2,...,v_N]^T$. Now  consider the central differences for $z''(x_j)\approx (-A_h\vec{v})_j$ with the matrix $A_h$ defined by
\begin{equation}\label{matrixa}
	A_h^{FD}:=\frac{c^2}{h^2}\begin{bmatrix}
		2&-1&0&\dots&\dots&\dots&0\\
		-1&2&-1&0&\dots&\dots&0\\
		&\ddots&\ddots&\ddots&\ddots&\ddots&\\
		0&\dots&\dots&0&-1&2&-1\\
		0&\dots&\dots&\dots&0&-1&1\\
	\end{bmatrix}_{N\times N}
	\end{equation}
	whose eigen-pairs  $(\mu^{FD}_{k}(h),\vec{\phi}_k(h))$ are \cite{B}:
	\begin{equation}\label{Inf}
		\left\{ \begin{array}{ll}
			\mu_k^{FD}(h)=\frac{4c^2}{h^2}\sin^2\left(\frac{(2k-1)\pi h}{2(2-h)L}\right),&\\
			\phi_{k,j}=\sin\left(\frac{(2k-1)j\pi h}{L(2-h)}\right),\quad k,j=1,2,...,N.&
		\end{array}\right.
	\end{equation}
	Considering no feedback $\xi\equiv 0$ in \eqref{estatic},  the following discretization  is obtained
\begin{equation}\label{FD}
	\left\{	\begin{array}{ll}
		\vec{\ddot{v}}+A_h^{FD}\vec{v}=0,\quad  t\in \mathbb{R}^+&\\
		v_0=0,\quad v_{N+1}-v_N=0,&\\
		v_j(0)=v^0_j,~\dot{v}_j(0)=v^1_j,\quad  j=0,...,N+1.
	\end{array}\right.
\end{equation}
The discretized  energy corresponding to \eqref{eq4} is
\begin{eqnarray}
\label{FD-en}
\begin{array}{ll}
	&E_{h,0}^{FD}(t):=\frac{h}{2}\sum\limits_{j=0}^{N}\rho\left|\dot v_j\right|^2+c^2\left|\frac{v_{j+1}-v_j}{h}\right|^2.
\end{array}
\end{eqnarray}
Defining   $\mathcal{A}_h^{FD}$ by $\mathcal{A}_h^{FD}(\vec{u}_1,\vec{u}_2):=\left(\vec{u}_2,-A_h^{FD}\vec{u}_1\right),$ and calling $\vec{y}_h=(\vec{v},\dot{\vec{v}})$,  \eqref{FD} can be reformulated as
\begin{equation}\label{semi-dis}
	\frac{d}{dt}\vec{y}_h=\mathcal{A}_h^{ FD}\vec{y}_h.
\end{equation}
\begin{Lem} \cite{B} For $K:=\{-N,...,-1,1,...N \},$ the eigen-pairs $\{(\lambda_k^{ FD}(h),\vec{\varphi}_k(h))\}_{K}$ of $\mathcal{A}_h^{ FD}$ are
	\begin{equation}\begin{array}{ll}
		(\lambda_k^{ FD},\vec{\varphi}_k)=\left(i\sqrt{\mu^{FD}_k(h)},\begin{bmatrix}
			\frac{1}{i\sqrt{\mu^{FD}_k(h)}}\vec{\phi}_k(h)\\
			\vec{\phi}_k(h)
		\end{bmatrix}\right),
		\end{array}
	\end{equation}
	and  $\sqrt{\lambda_k^{ FD}(h)}:=-\sqrt{\lambda_{-k}^{ FD}(h)}$ and $\vec{\varphi}_k:=\vec{\varphi}_{-k}$ for $k=-1,-2,...,-N.$ 	Therefore,  the solutions to \eqref{FD} can be expressed as
	\begin{equation}\label{discsol}
		\vec{v}(t)=\sum_{k\in K}(a_ke^{i\sqrt{\lambda_{k}^{ FD}(h)}t})\vec{\phi}_k(h).
	\end{equation}
\end{Lem}

 \subsection{Finite Elements with $\xi\equiv 0$}
First, multiply both sides of the equation in \eqref{estatic} by a continuously differentiable test function $\phi(x)\in C_0^\infty[0,L]$, and integrate both sides of the equation over $[0, L]$  to get
 \begin{eqnarray}\label{eq:1}
\int_{0}^{L} u_{tt} \phi~dx + \int_{0}^{L} u_{x} \phi_{x}~dx = 0. \\ \nonumber
\end{eqnarray}
At each node $\left\{x_i\right\}_{i=1}^N$, the following linear splines are defined 
\begin{eqnarray}
&&\phi_i(x)=\left\{\begin{array}{ll}
\frac{1}{h}(x-x_i),& x_{i-1}<x<x_{i}\\
\frac{1}{h}(x-x_{i+1}),& x_{i}<x<x_{i+1}\\
0,& {\textrm otherwise,}
\end{array}
\right.\\
&&\phi_{N+1}(x)=\left\{\begin{array}{ll}
\frac{1}{h}(x-x_N),& x_{N}<x<x_{N+1}\\
0,& {\textrm otherwise}.
\end{array}
\right.
\end{eqnarray}
Defining the $(N+1)\times (N+1)$ matrices $A_h$ and $M$ by
\begin{equation}\label{matrixa-fem}
	A_h^{FEM}:=\frac{c^2}{h^2}\begin{bmatrix}
		2&-1&0&\dots&\dots&\dots&0\\
		-1&2&-1&0&\dots&\dots&0\\
		&\ddots&\ddots&\ddots&\ddots&\ddots&\\
		0&\dots&\dots&0&-1&2&-1\\
		0&\dots&\dots&\dots&0&-1&1\\
	\end{bmatrix},
	\end{equation}
\begin{eqnarray}\label{mass} M:=\left(
                                                 \begin{array}{cccccc}
                                                   2/3 & 1/6 & 0 & 0& \ldots &  0\\
                                                   1/6 & 2/3 & 1/6 & 0&\ldots & 0 \\
                                                   0 & 1/6 & 2/3 & 1/6 &\ldots & 0\\
                                                   \vdots & \vdots & \vdots & \ddots &\vdots & \vdots\\
                                                   0 & 0 & 0 & 1/6 & 2/3 &1/6\\
                                                   0 & 0 & \ldots & 0 & 1/6 & 1/3\\
                                                 \end{array}
                                               \right),
\end{eqnarray}
and seeking solutions to  \eqref{eq:1} with $\xi\equiv 0$ of the form $v(x,t) = \sum\limits_{i=0}^{N+1} v_i(t)\phi_i(x)$
leads to
\begin{equation}\label{FEM}
	\left\{	\begin{array}{ll}
		\vec{\ddot{v}}+M^{-1}A^{FEM}_h\vec{v}=0,\quad  t\in \R^+&\\
		v_0=0,\quad h\frac{2\ddot v_{N+1}+\ddot v_N}{6}+c^2\frac{v_{N+1}-v_N}{h}=0,&\\
		v_j(0)=v^0_j,~\dot{v}_j(0)=v^1_j,\quad  j=0,...,N+1.
	\end{array}\right.
\end{equation}
The discretized energy corresponding to \eqref{FEM} is defined by
\begin{eqnarray}
\begin{array}{ll}\label{FEM-en}
	E_h^{ FEM}(t):=\frac{h}{12}\left[|\dot v_{N+1}|^2+ \sum_{j=1}^{N} \left(2\left|\dot v_j\right|^2+ \left|\dot v_j+\dot v_{j+1}\right|^2  + 6 c^2\left|\frac{v_{j+1}-v_j}{h}\right|^2\right)\right].
\end{array}
\end{eqnarray}
Defining   $\mathcal{A}_h^{ FEM}(\vec{u}_1,\vec{u}_2):=\left(\vec{u}_2,-M^{-1}A^{FEM}_h\vec{u}_1\right),$ and calling $\vec{y}_h=(\vec{v},\dot{\vec{v}})$,  \eqref{FEM} can be reformulated as
\begin{equation}\label{FEM2}
	\frac{d}{dt}\vec{y}_h=\mathcal{A}_h^{ FEM}\vec{y}_h.
\end{equation}
Introduce an $(N+1)\times (N+1)$ diagonal matrix $K$ by $K={\mathrm diag} (2,\ldots,2,1)$ so that \begin{eqnarray}\label{eigen}
\begin{array}{ll}
\lambda_j\left(M^{-1}A_h^{FEM}\right) = \lambda_j\left(M^{-1}KK^{-1}A_h^{FEM}\right) = \frac{ \lambda_j\left(K^{-1}A_h^{FEM}\right)}{ \lambda_j\left(K^{-1}M\right)}, ~j = 1,  \ldots, N+1,
\end{array}
\end{eqnarray}
where matrices $K^{-1}A_h$ and $K^{-1}M$ are band-matrices:
\begin{eqnarray}\label{mat0}
\begin{array}{ll}K^{-1}A_h^{FEM}=\frac{c^2}{2h^{2}}\left(
                                                 \begin{array}{cccccc}
                                                   2 & -1 & 0 & 0 & \ldots &  0\\
                                                   -1 & 2 & -1 & 0&\ldots & 0 \\
                                                   \vdots & \vdots & \vdots & \ddots &\vdots & \vdots\\
                                                   0 & 0 & 0 & -1 & 1 &-1\\
                                                   0 & 0 & \ldots & 0 & -2 & 2\\
                                                 \end{array}
                                               \right),\\
                                               K^{-1}M=\left(
                                        \begin{array}{cccccc}
                                        \frac{1}{3} & \frac{1}{12} & 0 & 0 & \ldots &  0\smallskip\\
                                        \frac{1}{12} & \frac{1}{3} & \frac{1}{12} & 0&\ldots & 0 \smallskip\\
                                        \vdots & \vdots & \vdots & \ddots &\vdots & \vdots\\
                                        0 & 0 & 0 & \frac{1}{12} & \frac{1}{3} &\frac{1}{12}\smallskip\\
                                        0 & 0 & \ldots & 0 & \frac{1}{6} & \frac{1}{3}\\
                                        \end{array}
                                    \right).
                                               \end{array}
\end{eqnarray}
The next two lemmas are necessary for finding the eigenvalues of matrices (\ref{mat0}).
\begin{Lem}\label{lem1}
The eigenvalues and eigenfunctions of the matrix $K^{-1}A_{h}^{FEM}$, respectively, are
	\begin{equation}\label{Inf}
		\left\{ \begin{array}{ll}
			\lambda_{j}(h)=\frac{2c^2}{h^2}\sin^2\left(\frac{(2j-1)\pi}{4N}\right),&\\
			u_{j,k}=\sin\left(\frac{(2j-1)k\pi}{2N}\right), ~ j,k=1,2,...,N+1.
		\end{array}\right.
	\end{equation}
\end{Lem}

{\bf Proof:} Letting $u=[u_{1},u_{2},\ldots,u_{N}]^{T}$, the eigenvalue problem for $k = 1, 2, \ldots, N+1$ is
\begin{eqnarray}\label{Inf-disc}
 \label{dif-eq}\left\{ \begin{array}{ll}
 -u_{k-1} + (2-2h^{2}\lambda)u_k - u_{k+1} = 0, &  \\
    u_0 = 0, ~~ u_{N+1} = u_{N}.
    \end{array}\right.
\end{eqnarray}
With $u_k = z^k$,  $u_{k-1} = z^{k-1}$ and $u_{k+1} = z^{k+1}$,
both $z^{k-1} \neq 0$ and $ (-1+(2-2h^{2}\lambda)z-z^{2})=0.$ Thus,
 the  general solution to \eqref{Inf-disc} is
\begin{eqnarray}\label{u-sub-k} & u_{k}=c_{1}z^{k}(\lambda)+c_{2}z^{-k}(\lambda).
\end{eqnarray}
By the boundary conditions, 
\begin{eqnarray}
\begin{array}{ll}
   & z_j = e^{\frac{i(2j-1)k\pi}{2N}}, ~j = 1, \ldots, N+1.\\
    \end{array}
\end{eqnarray}
Therefore, substituting $z_j$ into (\ref{u-sub-k})  leads to
\begin{eqnarray}
\begin{array}{ll}
u_{j,k} = sin\left(\frac{(2j-1)k\pi h}{2(L-h)}\right), \quad j,k = 1, 2, \ldots, N+1.
\end{array}
\end{eqnarray}
Lastly, the eigenvalues are solved by using $z_1z_2 = 1$ and $z_1 + z_2 = 2-2h^2\lambda_j$:
\begin{eqnarray}
2\left(1-2\sin^2\left({\frac{(2j-1)\pi}{4N}}\right)\right) = 2-2h^2\lambda_j.
\end{eqnarray}
Therefore, noting $N = \frac{L-h}{h}$, \eqref{Inf} is obtained. $\square$

\begin{Lem}\label{thm00000} For $ j = 1,  \ldots, N+1,$ the eigenvalues $M^{-1}A_h^{FEM}$ of \eqref{FEM2} are given by
\begin{eqnarray}\label{FEM_eigs}
{\lambda}_j(M^{-1} A_h^{FEM}) =  \frac{1}{h^2}\frac{6-6\cos{\left(\frac{(2j-1)\pi h}{2(L-h)}\right)}}{2+\cos{\left(\frac{(2j-1)\pi h}{2(L-h)}\right)}}.
\end{eqnarray}
\end{Lem}

{\bf Proof:} Define the $(N+1)\times (N+1)$ matrix $J$ by $J:={\mathrm {tridiag}(1,0,1)}.$ 
One can readily verify that $h^2K^{-1}A = I - K^{-1}J$  and $6K^{-1}M = 2I + K^{-1}J.$ Since $h^2K^{-1}A_h^{FEM}$ and $6K^{-1}M$ are diagonalizable, there exists an invertible matrix P such that $
6P^{-1}K^{-1}MP = 2I + P^{-1}K^{-1}JP.$
Lastly, letting $P^{-1}K^{-1}AP = D$ where D is a diagonal matrix,
$6P^{-1}K^{-1}MP = 3I - h^2D$, $P^{-1}K^{-1}MP = \frac{1}{6}\left(3I - h^2D\right).$
Therefore, the eigenvalues of the matrix $K^{-1}M$ lie on the diagonal of matrix $P^{-1}K^{-1}MP$,
\begin{eqnarray}
&\lambda_j\left(K^{-1}M\right) = \frac{1}{6}\left(3 - h^2\lambda_j\left(K^{-1}A_h^{FEM}\right)\right).~
\end{eqnarray}
Therefore,
\begin{eqnarray}
\begin{array}{ll}
\lambda_j\left(M^{-1}A_h^{FEM}\right)& = \frac{ \lambda_j\left(K^{-1}A_h^{FEM}\right)}{ \lambda_j\left(K^{-1}M\right)} 
= \frac{\frac{12}{h^2}\sin^2{\left(\frac{(2j-1)\pi h}{4(L-h)}\right)}}{2\left(\frac{3}{2}-\sin^2{\left(\frac{(2j-1)\pi h}{4(L-h)}\right)}\right)}.
\end{array}
\end{eqnarray}
Following the sub-eigenvalue problem (\ref{eigen}), \eqref{FEM_eigs} is obtained. $\square$

\begin{theorem} It can be shown that  the eigen-pairs $\{\lambda_k^{FEM}(h),\vec{\varphi}_k(h)\}_{K}$ of $\mathcal{A}_h$ are
	\begin{equation}
	\begin{array}{ll}
		\left(i\sqrt{\lambda_k(M^{-1}A_h^{FEM})},\begin{bmatrix}
			\frac{1}{i\sqrt{\lambda_k(M^{-1}A_h^{FEM})}}\vec{\phi}_k\\
			\vec{\phi}_k
		\end{bmatrix}\right),
		\end{array}
	\end{equation}
	where  $\sqrt{\lambda_k^{FEM}(h)}:=-\sqrt{\lambda_{-k}^{FEM}(h)}$ and $\vec{\varphi}_k:=\vec{\varphi}_{-k}$ for $k=-1,-2,...,-N.$ 	 The solutions to \eqref{FD} can be expressed as
	\begin{equation}\label{discsol}
		\vec{v}(t)=\sum_{k\in K}(a_ke^{i\sqrt{\lambda_{k}^{ FEM}(h)}t})\vec{\phi}_k(h).
	\end{equation}
\end{theorem}

	\section{Lack of Observability and Exponential Stability as $h\to 0$}
	It can be shown easily the observability result stated in Theorem \ref{observa} does not hold uniformly for the discretized models since
	\begin{eqnarray}
	\label{lackof}\left\{
	\begin{array}{ll}
	\frac{h^2|\lambda_{N+1}^{FD}(h)|^2}{c^2}\to 4,\qquad \frac{h^2|\lambda_{N+1}^{ FEM}(h)|^2}{c^2}\to 12 & {\rm as} ~N\to \infty.
	\end{array}\right.
	\end{eqnarray}
 \begin{Lem}\label{lem} Considering either  \eqref{FD} or \eqref{FEM}, for any $T>0,$ respectively,
 \begin{eqnarray}
\lim\limits_{h\to 0} \sup \frac{E_{h,0}^{FD}(0)}{\int_0^T \left|\frac{\dot v_{N}}{h}\right|^2},\qquad \lim\limits_{h\to 0} \sup \frac{E^{ FEM}_h(0)}{\int_0^T \left|\frac{\dot v_{N}}{h}\right|^2}\to \infty.
 \end{eqnarray}
		\end{Lem}
	{\bf Proof}: The first limit is proved in \cite{B}. The second limit follows from the same argument together with Theorem \ref{thm00000}. $\square$

	\subsection{Finite Differences with $\xi>0$}
First, observe that  the absorbing boundary conditions $c^2 v_x(L,t)=-\xi\dot v(L,t)$ in \eqref{estatic} can be approximated by
\begin{eqnarray*}
\begin{array}{ll}
c^2 v_{xx} (x_{N+1})& \approx c^2 \frac{v_x(x_{N+1})-v_x(x_N)}{h}=  \frac{-\xi \dot v_{N+1}-c^2 \left(\frac{v_{N+1}-v_{N}}{h}\right)}{h}.
\end{array}
\end{eqnarray*}
Then, the following Finite-Difference approximation for \eqref{estatic} can be considered
\begin{equation}\label{FD-exp}
	\left\{	\begin{array}{ll}
		\left\{\ddot v_j-c^2\left(\frac{v_{j+1}+v_{j-1}-2v_j}{h^2}\right)=0\right\}_{j=1}^N,\\
		v_0=0, \ddot v_{N+1}+c^2\frac{v_{N+1}-v_{N}}{h^2} +\xi \frac{\dot v_{N+1}}{h}=0,\\
		v_j(0)=v^0_j,~\dot{v}_j(0)=v^1_j,\quad  j=0,...,N+1.
	\end{array}\right.
\end{equation}
The energy $E_{h,0}^{FD}(t)$, defined in \eqref{FD-en} for the case $\xi\ne 0$, is now redefined as the following
\begin{eqnarray}
\label{FD-en-new}
\begin{array}{ll}
	&E_h^{FD}(t):=\frac{h}{2}\sum\limits_{j=0}^{N+1}\rho\left|\dot v_j\right|^2+c^2 \frac{h}{2}\sum\limits_{j=0}^{N}\left|\frac{v_{j+1}-v_j}{h}\right|^2.
\end{array}
\end{eqnarray}

Letting  $\vec{u}_{1,h}=(v_1,v_2,...,v_{N+1})^{\rm T}$, $\vec{u}_{2,h}=(\dot v_1,...,\dot v_{N+1})^{\rm T}, $
$\vec{y}_h=(\vec{u}_{1,h},\vec{u}_{2,h})^{\mathrm T},$ the system \eqref{FD-exp} is written in the first-order form
\begin{equation}\label{eqr2}
    {\dot {\vec{y}}_h}=\mathcal{A}_h^{FD}(\xi)\vec{y}_h=(\vec{y}_{2,h},-A_h\vec{y}_{1,h}+B_h\vec{y}_{2,h} )^{\mathrm T}
    \end{equation}
    where $A_h$ is the matrix defined by \eqref{matrixa-fem}, and$B_{N+1,N+1}=\frac{-\xi}{h}\ne 0$ and otherwise $B_{i,j}\equiv 0$ for any other $i,j.$

Consider the eigenvalue problem  for \eqref{eqr2}:
\begin{eqnarray} \label{eig-FD-f}
\mathcal{A}_h^{FD}(\xi)\vec{y}=\lambda^{FD}(\xi,h)\vec{y}. 
\end{eqnarray}
For $0<\xi<c,$ it can be shown that the real part of the eigenvalues are negative, e.g. ${\rm Re} \lambda^{FD}(\xi,h)=\frac{-L\xi|\lambda^{FD} (\xi,h) y_{1,N+1}|^2}{hc \left(|\lambda^{FD}(\xi,h) \vec y_{1,h}|^2-\vec y_{1,h}^{\rm T} A \vec y_{1,h})\right)}<0.$

Seeking the solution of \eqref{eig-FD-f}  of the form $y_{1,k}=z^{2k}-z^{-2k}$ with $k=1,2,\ldots, N+1$ and $z\in \mathbb{C}$ leads to $\lambda^{FD}(\xi,h)=c(N+1)(z-z^{-1})$ where $z$ satisfies  $z^2\ne \mp 1$ and
\begin{equation}\label{la-pol}
   p(z)=z^{4N+6}+\frac{\xi}{c}z^{4N+5}-\frac{\xi}{c}z+1=0.
\end{equation}

Note that the eigenvalues $\lambda^{FD}(\xi,h)$ of  $\mathcal{A}_h^{FD}(\xi)$ solve \eqref{la-pol}, but some of the solutions of \eqref{la-pol} are not related to the eigenvalues  $\lambda^{FD}(\xi,h)$. In fact, since   $\mathcal{A}_h^{FD}(\xi)$ is a real matrix, it is sufficient to identify all eigenvalues of  $\mathcal{A}_h^{FD}(\xi)$ with positive imaginary part.

Moreover, observe that $z=\mp i$ are the roots of $p.$ On the other hand, if $z$ is a root of $p(z),$ then $-z$ and $z^{-1}$ are both the roots of $p(z).$ As $z\ne -\frac{\xi}{c},$ \eqref{la-pol} is equivalent to 
\begin{equation}\label{la-pol1}
   z^{4N+5}=\frac{\xi z-c}{\xi+c z}.
\end{equation}
The following list of results are proved in \cite{Peichl}.
\begin{Lem} \cite[Prop. 3.3]{Peichl}
The polynomial $p(z)$ satisfies the following:
\begin{itemize}
\item [1.] The only roots of $p(z)$ on the unit circle are $z=\mp i.$ Moreover, $z=\mp i$ are the only root of $p(z)$ with purely imaginary parts.
\item [2.]  All roots of $p(z)$ with positive real parts must be within the unit circle. All roots of $p(z)$ with negative  real parts must be outside the unit circle.\end{itemize}
\end{Lem}
It is inferred from these results that each root $z$ of $p(z)$ in the first quadrant of the complex plane leads to three other roots of $p(z): $ $\bar {z}, z^{-1}, {\bar z}^{-1}.$ Therefore, the eigenvalue of  $\mathcal{A}_h^{FD}(\xi)$  corresponding to the roots of $p(z)$ has positive imaginary part if and only if $z$ has positive imaginary part. As a result, we only estimate the roots of $p(z)$ in the first quadrant. 	

Following the clever discussion in \cite{Peichl},  the analysis of the distribution of the roots of \eqref{la-pol} is solely based on the analysis of the mapping defined by $\mathcal T(z):={\mathcal G}^{\frac{1}{4N+5}}(z)$ where ${\mathcal G}(z)=\frac{\xi z-c}{zc+\xi}.$ Indeed,  a fixed point $z$ of $\mathcal T$ is a root of the polynomial \eqref{la-pol}. Therefore, the root of $p$ satisfies $$z=[\mathcal{G}(z)]^{\frac{1}{4N+5}}=\left(\frac{\xi z-c}{zc+\xi}\right)^{\frac{1}{4N+5}}.$$

Define the following sector in the first quadrant of $\mathbb C:$
$$S=\{z \in \mathbb{C}~|~Re(z)\geq 0, Im(z)\geq 0, |z|\leq \frac{\xi+c}{2\xi}\}.$$
Since $T$ is multi-valued, considering the branches ${\mathcal T}_j$ of $\mathcal T$ as the following
\begin{eqnarray}
\begin{array}{ll}
&{\mathcal T}_j(z)=\mathcal{G}(z)]^{\frac{1}{4N+5}} e^{\frac{i(\theta(z)+2j\pi )}{4N+5}}, ~~1<j\le 4N+4,\\
&\theta(z)={\textrm Arg}({\mathcal G}(z)).
\end{array}
\end{eqnarray}
\begin{Lem} \cite[Prop. 3.4 \& 3.5, Cor. 3.1]{Peichl}
The following results hold:
\begin{itemize}
\item [1.] The complex functions ${\mathcal T}_j,$ $j=1,\ldots 4N+4,$ are analytic over $S,$ and  ${\mathcal T}_j$ is a contraction mapping over $S$ for $N$ large enough.
\item [2.] Since $G$ is analytic over $S,$ there exists a constant $M_{G}>0$ such that $|G(z)|<M_G$ for all $z\in S.$ 
\item [3.]  As well, for every $j\in [0, N],$ the subsections $S_j$ of $S$ defined by
$$S_j:=\left\{z\in S, {\textrm Arg} (z)\in \left[\frac{2j\pi}{4N+5},\frac{(2j+1)\pi }{4N+5}\right]\right\}$$
is invariant under ${\mathcal T}_j$ for  large enough $N$. 
\item [4.]
The fixed point $z$ of $\mathcal T_j$ in $S_j$ satisfies\\
\begin{eqnarray}
\label{zestimates1}  &&|z|\geq(\frac{\xi c-\xi^2}{2\xi^2+\xi c+c^2})^{\frac{1}{4N+5}}-\sqrt{2}(1+\frac{c}{\xi})\frac{{M_G}^{\frac{1}{4N+5}}}{4N+5},\\
    \label{zestimates2}&&|z|\leq \left(\frac{c}{\xi}\right)^{\frac{1}{4N+5}}+\left(1+\frac{c}{\xi}\right)\frac{{M_G}^{\frac{1}{8N+10}}}{4N+5}.
    \end{eqnarray}
    \end{itemize}
\end{Lem}
Now, notice that
\begin{eqnarray}
\begin{array}{ll}
{\mathrm Re} \lambda_j^{ FD}(\xi,h)=\frac{c(N+1)}{L} (z_j-\frac{1}{z_j}) \cos{({\textrm Arg} ~z_j)},\\
{\mathrm Im} \lambda_j^{ FD}(\xi,h)=\frac{c(N+1)}{L} (z_j+\frac{1}{z_j}) \sin{({\textrm Arg} ~z_j)},
\end{array}
\end{eqnarray}
where $Arg (z_j)=Arg(T_j(z_k))=\frac{\theta(z_j)+2j\pi}{4N+5}, j=0,1,\ldots 4N+4,$ and 
\begin{equation}\label{thetaz}
\begin{array}{ll}
   & \theta(z_j)=\pi-\arctan\left(\frac{Im ~\mathcal{G}(z_j)}{Re ~\mathcal{G}(z_j)}\right),\qquad \frac{\pi}{2}\le \theta(z_j)\le  \frac{3\pi}{2}.
    \end{array}
    \end{equation}
Finally, the following results describe the eigenvalues of $\mathcal{A}_h^{FD}(\xi)$ in terms of the roots of \eqref{la-pol}.
\begin{theorem}  \cite[Thm 3.1]{Peichl} For sufficiently large $N$, all roots of $p(z)$ are simple. There are exactly $N+1$ roots in the first quadrant of the complex plane. Each subset $S_k$ contains exactly one root $z_j$ of $p(z)$ for $j=1,2,\ldots N.$ As a result, the matrix $\mathcal{A}_h^{FD}(\xi)$ has $2N+2$ eigenvalues $\lambda_j^{FD}(\xi,h).$ 
Moreover, this implies that $Re ~\lambda_j^{FD}(\xi,h)\rightarrow0$ as $N\to\infty$. Hence, the system \eqref{FD-exp} lacks of exponential stability uniformly as $N\to \infty$. 
 \end{theorem}

Now we are ready to state one of our main results, which helps estimate the eigenvalues further:
 
  \begin{theorem}\label{thm-eigFD}
 For all $j=1,2,\ldots, 2N+2$ and   sufficiently large  $N$ (or small enough $h$), the eigenvalues $\lambda_j^{FD}(\xi,h)$ of $\mathcal{A}_h^{FD}(\xi)$ satisfy 
 \begin{eqnarray}\label{thm5.1}
    && h ~|Re \lambda_j^{FD}(\xi,h)|  =O(h),\\
   \label{thm5.2}&&h ~|Im \lambda_j^{FD}(\xi,h)|  =2c \left|\sin\left(\frac{2j\pi}{4N+5}\right)\right|+O(h). 
 \end{eqnarray}
 Hence, 
  \begin{equation}\label{impw} h^2| \lambda_j^{FD}(\xi,h)|^2\le 4c^2  \sin^2 \left( \frac{2j\pi}{4N+5}\right)+ O(h).
  \end{equation}
 \end{theorem}

{\bf Proof:} First observe that \eqref{thm5.1} holds by \eqref{thetaz}  since
\begin{eqnarray*}
    \begin{array}{ll}
       |h ~Re \lambda_j^{FD}(\xi,h)|&=ch\frac{N+1}{L}|(|z_j|-|z_j|^{-1})\cos(Arg(z_j))| \\
         & =ch\frac{1}{h}\left|\left(\frac{|z_j|^2-1}{|z_j|}\right)\cos \left(\frac{2j\pi}{4N+5}+\frac{\theta(z_j)}{4N+5}\right)\right|\\
         &=c\left|\frac{-1-|z_j|}{|z_j|}\right||1-|z_j||\left|\cos\left(\frac{2j\pi}{4N+5}\right)\cos\left(\frac{\theta(z_j)}{4N+5}\right)-\sin\left(\frac{2j\pi}{4N+5}\right)\sin\left(\frac{\theta(z_j)}{4N+5}\right)\right|,
    \end{array}
\end{eqnarray*}
and by utilizing \eqref{zestimates1} and \eqref{zestimates2} leading to \eqref{thm5.1}
\begin{eqnarray*}
    \begin{array}{ll}
      |h ~Re \lambda_j^{FD}(\xi,h)|&=ch\frac{N+1}{L}|(|z_j|-|z_j|^{-1})\cos(Arg(z_j))| \\
         & =ch\frac{1}{h}\left|\left(\frac{|z_j|^2-1}{|z_j|}\right)\cos \left(\frac{2j\pi}{4N+5}+\frac{\theta(z_j)}{4N+5}\right)\right|\\
         &=c\left|\frac{-1-|z_j|}{|z_j|}\right||1-|z_j||\left|\cos\left(\frac{2j\pi}{4N+5}\right)\cos\left(\frac{\theta(z_j)}{4N+5}\right)-\sin\left(\frac{2j\pi}{4N+5}\right)\sin\left(\frac{\theta(z_j)}{4N+5}\right)\right|\\
         &\leq c\left|-1-\left(\frac{\xi c-\xi^2}{2\xi^2+\xi c+c^2}\right)^{\frac{1}{4N+5}}+\sqrt{2}(1+\frac{c}{\xi})\frac{M^{\frac{1}{4N+5}}}{4N+5}\right|\\
         &\quad \times\left|1-\left(\frac{\xi c-\xi^2}{2\xi^2+\xi c+c^2}\right)^{\frac{1}{4N+5}}+\sqrt{2}(1+\frac{c}{\xi})\frac{M^{\frac{1}{4N+5}}}{4N+5}\right| \times|1+O(h)|\\
         &=c\left|-1-1+O(h)\right|O(h)\\
         &=O(h).
    \end{array}
\end{eqnarray*}

On the other hand, analogously, by \eqref{zestimates1}-\eqref{thetaz}
\begin{eqnarray*}
    \begin{array}{ll}
       \left|h~Im \lambda_j^{FD}(\xi,h)\right|  &= h \frac{c(N+1)}{L}|(|z_j|+|z_j|^{-1})\sin(Arg(z_j))|  \\
         &= c||z_j|+|z_j|^{-1}| \left|\sin \left(\frac{2j\pi}{4N+5}+\frac{\theta(z_j)}{4N+5}\right)\right|\\
         &= c||z_j|+|z_j|^{-1}|\left|\sin\left(\frac{2j\pi}{4N+5}\right)\cos\left(\frac{\theta(z_j)}{4N+5}\right)+\cos\left(\frac{2j\pi}{4N+5}\right)\sin\left(\frac{\theta(z_j)}{4N+5}\right)\right|\\
         &= c||z_j|+|z_j|^{-1}|\left|\sin\left(\frac{2j\pi}{4N+5}\right)+O(h)\right|\\
         &\leq c\left|\left(\frac{c}{\xi}\right)^{\frac{1}{4N+5}}+\left(1+\frac{c}{\xi}\right)\frac{M^{\frac{1}{8N+10}}}{4N+5}+\frac{1}{\left(\frac{\xi c-\xi^2}{2\xi^2+\xi c+c^2}\right)^{\frac{1}{4N+5}}-\sqrt{2}(1+\frac{c}{\xi})\frac{M^{\frac{1}{4N+5}}}{4N+5}}\right|\\
         &\qquad \qquad\times \left|\sin\left(\frac{2j\pi}{4N+5}\right)+O(h)\right|\\
         &\leq c\left|\sin\left(\frac{2j\pi}{4N+5}\right)+O(h)\right|\left|\left(\frac{c}{\xi}\right)^{\frac{1}{4N+5}}+O(h)+\frac{1}{\left(\frac{\xi c-\xi^2}{2\xi^2+\xi c+c^2}\right)^{\frac{1}{4N+5}}+O(h)}\right|\\
         &\leq c\left|\sin\left(\frac{2j\pi}{4N+5}\right)+O(h)\right|\left|1+O(h)+\frac{1}{1+O(h)}\right|\\
        &\le 2c \left|\sin\left(\frac{2j\pi}{4N+5}\right)\right|+O(h).
    \end{array}
\end{eqnarray*}
Finally, \eqref{impw} follows from \eqref{thm5.1} and \eqref{thm5.2}. $\square$

\subsection{Finite Elements with $\xi>0 $}

Considering the discretized model  by Finite Elements
\begin{equation}\label{FEM-exp}
	\left\{	\begin{array}{ll}
		\left\{\frac{\ddot v_{j+1} + 4 \ddot v_j+\ddot v_{j-1}}{6}+c^2\frac{v_{j+1} -2 v_j+ v_{j-1}}{h^2}=0\right\}_{j=1}^N,  &\\
		v_0=0,~\frac{2\ddot v_{N+1}+\ddot v_N}{6}+c^2\frac{v_{N+1}-v_N}{h^2}=-\frac{\xi}{h} \dot v_{N+1},&\\
		v_j(0)=v^0_j,~\dot{v}_j(0)=v^1_j,\quad  j=0,...,N+1,
	\end{array}\right.
\end{equation}
the analysis and arguments above can be replicated analogously, see \cite{Peichl} for more in-depth discussion, to get the following result.
 \begin{theorem}\label{thm-eigFEM}
 For all $j=1,2,\ldots, 2N+2$ and  sufficiently  large enough  $N$ (or small enough $h$ ), the eigenvalues $\lambda_j^{FEM}(\xi,h)$ of $\mathcal{A}_h^{FEM}(\xi)$ satisfy 
  \begin{equation}\label{imp1} h^2| \lambda_j^{FEM}(\xi,h)|^2\le 12c^2  \sin^2 \left( \frac{2j\pi}{4N+5}\right)+ O(h).
  \end{equation}
 \end{theorem}

\section{Exponential Stability  as $h\to 0$}

\subsection{Finite Differences with $\xi\ne 0$}
\begin{Lem} The system \eqref{FD-exp} is dissipative, i.e.
\begin{eqnarray}\label{FD-exp5}
\frac{dE_h^{ FD}}{dt} +\xi |\dot v_{N+1}|^2=0.\end{eqnarray}
\end{Lem}
{\bf Proof:} Multiply both sides of \eqref{FD-exp} by $h\dot v_j$ and take sum from $j=1$ to $N:$
\begin{eqnarray}\label{FD-exp2}
		\begin{array}{ll}
		\sum\limits_{j=1}^N h\dot v_j\ddot v_j -c^2\sum\limits_{j=1}^N \frac{v_{j+1} -2 v_j+ v_{j-1}}{h}\dot v_j=0.
	\end{array}&&
\end{eqnarray}
Since 
    \begin{eqnarray}
\label{dum}\begin{array}{ll}
	&-h c^2\sum\limits_{j=1}^N \frac{v_{j+1} -2 v_j+ v_{j-1}}{h^2}\dot v_j=\frac{c^2}{h} (v_N-v_{N+1}) \dot v_{N+1}  +\frac{c^2}{h} \sum\limits_{j=0}^N(v_{j+1}-v_{j}) (\dot v_{j+1}-\dot v_{j}),
\end{array}
\end{eqnarray}
substituting \eqref{dum} into \eqref{FD-exp2} yields 
\begin{eqnarray}\label{FD-exp222}
		\begin{array}{ll}
		\frac{dE_h^{ FD}}{dt}  -\frac{c^2}{h} (v_{N+1}-v_{N}) \dot v_{N+1}-h\dot v_{N+1} \ddot v_{N+1}=0,	\end{array}&&
\end{eqnarray}
and this together with the boundary conditions \eqref{FD-exp},
 \eqref{FD-exp5} follows. $\square$

Following \eqref{dm1},  define the following Lyapunov functional
\begin{eqnarray}\label{FD-exp6}
L_h^{FD}(t):=E_h^{FD}+ \delta F_h^{FD}(t)
\end{eqnarray}
where the auxiliary function $F_h^{ FD}$ is defined by
\begin{eqnarray}\label{imp10}
\begin{array}{ll}
    F_h^{FD}(t)=\sum\limits_{j=1}^N jh\dot v_j\left(\frac{v_{j+1}-v_{j-1}}{2}\right)+\frac{L}{2}(v_{N+1}-v_N)\dot v_{N+1} -\frac{L\xi h}{4c^2} |\dot v_{N+1}|^2.
    \end{array}
\end{eqnarray}
\begin{Lem} \label{lem 2} For $0<\delta <\frac{c}{L},$ $L_h^{FD}(t)$ is equivalent to $E_h^{FD}(t),$
\begin{equation}\label{eqrfd2}
\left(1-\frac{L\delta}{c}\right)E_h^{FD}\le L_h^{ FD}\le \left(1+\frac{L\delta}{c}\right) E_h^{ FD}.
\end{equation}
\end{Lem}

{\bf Proof:} First of all, by H\"olders, and Cauchy-Schwartz inequalities for sums,  the first and second terms in \eqref{imp10} are estimated as the following
\begin{eqnarray}
\label{eqrfd221}
\begin{array}{ll}
  &  \left|\sum\limits_{j=1}^N jh\dot v_j\left(\frac{v_{j+1}-v_{j-1}}{2}\right)\right| \le \frac{Lh}{2c}\left[ \sum\limits_{j=0}^{N+1} |\dot v_{j+1}|^2+\sum\limits_{j=0}^Nc^2\left |\frac{v_{j+1}-v_{j}}{h}\right|^2-\frac{c^2}{2}\left |\frac{v_{N+1}-v_{N}}{h}\right|^2-|\dot v_{N+1}|^2\right],
 \end{array}\\
\label{eqrfd222}
\begin{array}{ll}
&  \frac{L}{2} \left|(v_{N+1}-v_N)\dot v_{N+1} \right|\le  \frac{Lh}{4c} \left[c^2\left |\frac{v_{N+1}-v_{N}}{h}\right|^2 +\dot v_{N+1}^2\right].
 \end{array}
\end{eqnarray}
Considering \eqref{eqrfd221} and \eqref{eqrfd222}, the following is immediate
\begin{eqnarray}
\label{eqrfd22}
\begin{array}{ll}
    |F_h^{FD}|&
     \leq  \frac{Lh}{c} E_h(t) -\frac{Lh}{4c}|\dot v_{N+1}|^2-\frac{L\xi h}{4c^2} |\dot v_{N+1}|^2\leq   \frac{Lh}{c} E_h(t).
\end{array}
\end{eqnarray}
Therefore, \eqref{eqrfd2} follows from \eqref{eqrfd22}. $\square$
\begin{Lem}  The function $F_h$ satisfies
\begin{eqnarray}
\label{FD-exp7}
\begin{array}{ll}
&\frac{dF_h^{FD}(t)}{dt}  \le - \left(1-\frac{\kappa^{FD} }{4c^2}\right)E_h^{FD}(t) +\frac{L}{2}\left( 1+\frac{1}{2N+2} +\frac{\xi^2}{c^2} \right)|\dot v_{N+1}|^2,
 \end{array}
 \end{eqnarray}
where $\kappa^{FD} :=\max \left\{h^2\left|\lambda_j^{FD}(\xi, h)\right|^2\right\}_{j=1}^{N+1},$ and the upper bound for the term $h^2\left|\lambda_j^{FD}(\xi, h)\right|$ is given by  \eqref{impw}.
\end{Lem}
{\bf Proof}:  Finding the derivative of $F_h(t)$ along the solutions of \eqref{FD-exp} leads to
\begin{eqnarray}
 \begin{array}{ll}
    \frac{dF_h^{FD}}{dt}  &=   h\left[\sum\limits_{j=1}^N j\ddot{v}_j (\frac{v_{j+1}-v_{j-1}}{2})+\dot v_j\frac{{\dot v}_{j+1}- {\dot v}_{j-1}}{2})\right]+\frac{L}{2}(\dot v_{N+1}-\dot v_N)(\dot v_{N+1})\\
    & \quad +\frac{L}{2}( v_{N+1}- v_N)(\ddot v_{N+1})-\frac{L\xi h}{2c^2} \dot v_{N+1} \ddot v_{N+1}\\
    &= -\frac{hc^2}{2}\sum\limits_{j=0}^N |\frac{v_{j+1}-v_j}{h}|^2-\frac{h}{2}\sum\limits_{j=0}^{N}|{\dot v}_j|^2+\frac{h^3}{4} \sum\limits_{j=0}^N |\frac{\dot v_{j+1}-\dot v_j}{h}|^2 +\frac{Lc^2}{2}\left|\frac{v_{N+1}- v_N}{h}\right|^2-\frac{h}{4}|{\dot v}_{N+1}|^2\\
    &\quad -\frac{Lh^2}{4}\frac{|{\dot v}_{N+1}-{\dot v}_N|^2}{h^2}  + \frac{L(|{\dot v}_{N+1}|^2+|{\dot v}_N|^2)}{4}  +\frac{L}{2}(\dot v_{N+1}-\dot v_N)(\dot v_{N+1}) \\
    &\quad +\frac{L}{2} (v_{N+1}-v_N) \ddot v_{N+1} -\frac{L\xi h}{2c^2} \dot v_{N+1} \ddot v_{N+1}\\
        &= -E_h^{FD}(t)+\frac{h}{2}|\dot v_{N+1}|^2+\frac{h^3}{4} \sum\limits_{j=0}^N |\frac{\dot v_{j+1}-\dot v_j}{h}|^2 -\frac{h}{4}|{\dot v}_{N+1}|^2\\
    &\quad  +\frac{L}{2}|\dot v_{N+1}|^2-\frac{L\xi}{2h } (v_{N+1}- v_N)\dot v_{N+1}  -\frac{L\xi h}{2c^2} \dot v_{N+1} \ddot v_{N+1}\\
              &= -E_h^{FD}(t)+\frac{h^3}{4} \sum\limits_{j=0}^N |\frac{\dot v_{j+1}-\dot v_j}{h}|^2+\frac{h}{4}|{\dot v}_{N+1}|^2   -\frac{L\xi h}{2c^2} \dot v_{N+1} \ddot v_{N+1}\\
    &\quad+\frac{L}{2} \dot v_{N+1} \left[\dot v_{N+1}+\frac{\xi h}{c^2}\ddot v_{N+1}+\frac{\xi^2}{c^2} \dot v_{N+1}\right] \\
    &= -E_h^{FD}(t)+\frac{h^3}{4} \sum\limits_{j=0}^N |\frac{\dot v_{j+1}-\dot v_j}{h}|^2 +\frac{L}{2}\left(1+\frac{1}{2N+2}+\frac{\xi^2}{c^2}\right) |\dot v_{N+1}|^2.
     \end{array}
\end{eqnarray}
Defining the higher-order energy $\tilde E^{FD}_h(t):=\frac{h}{2}\sum\limits_{j=0}^{N}\rho\left|\ddot v_j\right|^2+c^2\left|\frac{\dot v_{j+1}-\dot v_j}{h}\right|^2
,$ the inequality above is reduced to 
\begin{eqnarray}\nonumber
 \begin{array}{ll}
    \frac{dF_h^{FD}}{dt} &\le  -E_h^{FD}(t)+\frac{h^2}{4c^2} {\tilde E}_h^{FD}(t) +\frac{L}{2}\left(1+\frac{1}{2N+2}+\frac{\xi^2}{c^2}\right) |\dot v_{N+1}|^2
     \end{array}
\end{eqnarray}
It follows from the definition of $\kappa_{FD}$ that $h^2\tilde E^{FD}_h(t)\le \kappa^{FD} E^{FD}_h(t)$. Hence, \eqref{FD-exp7} follows. $\square$ 

\begin{theorem} \label{main1}Suppose that there exists a constant $0<\delta <\frac{2\xi c^2}{L\left[c^2\left(1+\frac{1}{2N+2}\right)+ \xi^2\right]}<\frac{2\xi c^2}{L\left[c^2+ \xi^2\right]}$ such that for all initial conditions $\vec v^0, \vec v^1\in \mathbb{R}^{N+2},$ the energy $E_h^{ FD}(t)$ corresponding to \eqref{FD-exp}  satisfies
\begin{eqnarray}\label{thm-FD}
E_h^{FD}(t)\le \frac{c+\delta L}{c-\delta L} e^{-\delta(1-\frac{L\delta}{c})(1-\frac{\kappa^{FD}}{4c^2})t} E_h^{FD}(0), \quad \forall t>0.
\end{eqnarray}
\end{theorem}
{\bf Proof}: Since $\frac{L_h^{FD}(t)}{dt}=\frac{dE^{FD}(t)}{dt}+\delta \frac{dF_h^{FD}(t)}{dt},$ by Lemma \ref{lem 2}
\begin{eqnarray}
\nonumber \begin{array}{ll}
\frac{dL_h^{FD}(t)}{dt} &=- \delta\left(1-\frac{\kappa^{FD} }{4c^2}\right)E_h^{FD}(t) -\left(\xi - \frac{\delta L}{2}\left[1+ \frac{1}{2N+2}+ \frac{\xi^2}{c^2}\right]\right) |\dot v_{N+1}|^2\\
 &\le - \delta \left(1-\frac{L\delta}{c}\right)\left(1-\frac{\kappa^{FD}}{4c^2}\right)L_h^{ FD}(t).
 \end{array}
\end{eqnarray}
By the Gr\"onwall's inequality,
\begin{eqnarray}\label{imp1}
\begin{array}{ll}
L_h^{FD}(t)\le e^{-  \delta \left(1-\frac{L\delta}{c}\right)\left(1-\frac{\kappa^{FD}}{4c^2}\right)t}L_h^{FD}(0),
 \end{array}
\end{eqnarray}
and Lemma \ref{lem 2}, \eqref{thm-FD} is obtained. $\square$

\subsection{Finite Elements with $\xi\ne 0$}

Letting  $\vec{u}_{1,h}=(v_1,v_2,...,v_{N+1})^{\rm T}$, $\vec{u}_{2,h}=(\dot v_1,...,\dot v_{N+1})^{\rm T}, $
$\vec{y}_h=(\vec{u}_{1,h},\vec{u}_{2,h})^{\mathrm T},$ the system \eqref{FEM-exp} is written in the first-order form
\begin{equation}\label{eqr22}
    {\dot {\vec{y}}_h}=\mathcal{A}_h^{FEM}(\xi)\vec{y}_h=(\vec{y}_{2,h},-M^{-1}A_h^{FEM}\vec{y}_{1,h}+B_h\vec{y}_{2,h} )^{\mathrm T}
    \end{equation}
    where $A_h$ is  defined by \eqref{matrixa}, and $B_{N+1,N+1}=\frac{-\xi}{h}\ne 0$ and otherwise $B_{i,j}\equiv 0$ for any other $i,j.$

Consider the eigenvalue problem  for \eqref{eqr2}:
\begin{eqnarray} \label{eig-FEM-f}
\mathcal{A}_h^{FEM}(\xi)\vec{y}=\lambda^{FEM}(\xi,h)\vec{y}. 
\end{eqnarray}
For $0<\xi<c,$ it can be shown that the real part of the eigenvalues are negative, e.g. ${\rm Re} \lambda^{FEM}(\xi,h)=\frac{-L\xi|\lambda^{FEM} (\xi,h) y_{1,N+1}|^2}{hc \left(|\lambda^{FEM}(\xi,h) \vec y_{1,h}|^2-\vec y_{1,h}^{\rm T} A \vec y_{1,h})\right)}<0.$

\begin{Lem} The system \eqref{FEM}, or \eqref{eqr2}, is dissipative, i.e.
\begin{eqnarray}\label{FEM-exp5}
\frac{dE_h^{ FEM}}{dt} +\xi |\dot v_{N+1}|^2=0.\end{eqnarray}
\end{Lem}
{\bf Proof}: Multiply both sides of \eqref{FEM-exp} by $h\dot v_j$ and take sum from $j=1$ to $N:$
\begin{eqnarray}\label{FEM-exp2}
		\begin{array}{ll}
		\sum\limits_{j=1}^N\frac{\ddot v_{j+1} + 4 \ddot v_j+\ddot v_{j-1}}{6}h\dot v_j +c^2\sum\limits_{j=1}^N \frac{v_{j+1} -2 v_j+ v_{j-1}}{h}\dot v_j=0,
	\end{array}&&
\end{eqnarray}
where
\begin{eqnarray}
\label{FEM-exp3} \begin{array}{ll}
	&	h\sum\limits_{j=1}^N\frac{\ddot v_{j+1} + 4 \ddot v_j+\ddot v_{j-1}}{6}\dot v_j=-\frac{h}{6} (\ddot v_{N+1}+\ddot v_{N})\dot v_{N+1} +\frac{h}{3}\sum\limits_{j=0}^N  \ddot v_j \dot v_j+\frac{h}{6} \sum\limits_{j=0}^N (\dot v_{j+1} +\dot v_j) 	(\ddot v_{j+1}+\ddot v_j),
\end{array} &&\\
\nonumber  &&\\
\label{FEM-exp4}
\begin{array}{ll}
	&-h c^2\sum\limits_{j=1}^N \frac{v_{j+1} -2 v_j+ v_{j-1}}{h^2}\dot v_j=\frac{c^2}{h} (v_N-v_{N+1}) \dot v_{N+1}   +\frac{c^2}{h} \sum\limits_{j=0}^N(v_{j+1}-v_{j}) (\dot v_{v+1}-\dot v_{j}).
\end{array}
\end{eqnarray}
Substituting \eqref{FEM-exp3} and \eqref{FEM-exp4} into \eqref{FEM-exp2} yields  \eqref{FEM-exp5}. $\square$

Following \eqref{dm1}, define the following Lyapunov functional
\begin{eqnarray}\label{FEM-exp6}
L_h^{ FEM}(t):=E_h^{ FEM}+ \delta F_h^{ FEM}(t)
\end{eqnarray}
where the auxiliary function $F_h(t)$ is defined by
\begin{eqnarray}
\begin{array}{ll}
F_h^{ FEM}(t)=h\sum\limits_{j=1}^N\frac{\dot v_{j+1} + 4 \dot v_j+\dot v_{j-1}}{6} j\frac{ v_{j+1}-v_{j-1}}{2} +\frac{L}{6}\left(2\dot v_{N+1}+\dot v_{N}\right)(v_{N+1}-v_N).
 \end{array}
\end{eqnarray}
\begin{Lem} \label{lem2}For $0<\delta <\frac{c}{L},$ $L_h$ is equivalent to $E_h^{ FEM}, $ i.e.
\begin{equation}\label{}
\left(1-\frac{L\delta}{c}\right)E_h^{ FEM}\le L_h^{ FEM}\le \left(1+\frac{L\delta}{c}\right) E_h^{ FEM}.
\end{equation}
\end{Lem}
{\bf Proof:}
Applying the Cauchy-Schwartz inequality leads to
\begin{eqnarray}
\begin{array}{ll}
   |F_h^{ FEM}(t)|&\le \frac{L h}{6} \left|\frac{2\dot v_{N+1}+\dot v_N}{2}\right|^2+\frac{L h}{6} \left|\frac{v_{N+1}- v_N}{h}\right|^2\\
   & ~~+\frac{L h}{2} \sum\limits_{j=1}^N \left|\frac{v_{j+1}-v_j}{2h}\right|^2+\frac{Lh}{2} \sum\limits_{j=1}^N \left|\frac{\dot v_{j+1}+4\dot v_j+\dot v_{j-1}}{6}\right|^2\\
 &\le \frac{Lh}{9} \left|\frac{\dot v_{N+1}+\dot v_N}{2}\right|^2+\frac{Lh}{4} \left|\frac{v_{N+1}- v_N}{h}\right|^2+\frac{Lh}{18} |\dot v_{N+1}|^2+\frac{Lh}{4} \sum\limits_{j=1}^N \left|\frac{v_{j+1}-v_j}{h}\right|^2\\
 & ~~+\frac{Lh}{4} \sum\limits_{j=1}^N \left|\frac{v_j-v_{j-1}}{h}\right|^2+\frac{Lh}{9} \sum\limits_{j=1}^N \left|\frac{\dot v_{j+1}+\dot v_j}{2}\right|^2+\left|\frac{\dot v_{j}+\dot v_{j-1}}{2}\right|^2+\left|\dot v_j\right|^2\\
 & \le\frac{L}{c}E_h^{ FEM}(t).
\end{array}
\end{eqnarray}
Therefore,  \eqref{FEM-exp6} follows. $\square$

\begin{Lem} The function $F_h^{ FEM}$ satisfies
\begin{equation}\label{FEM-exp7}
\begin{array}{ll}
 \frac{dF_h^{ FEM}}{dt}\le& - \left(1-\frac{\kappa^{ FEM} }{12}\right)E_h^{ FEM}(t)+ \frac{L}{2}\left(1+\frac{\xi^2}{c^2}\right) |\dot v_{N+1}|^2
 \end{array}
\end{equation}
where $\kappa_{ FEM} =\max \left\{h^2\left|\lambda_j^{ FEM}(\xi,h)\right|^2\right\}_{j=1}^{N+1}.$
\end{Lem}
{\bf Proof}: Finding the derivative of $F^{ FEM}_h(t)$ along the solutions of \eqref{FEM-exp} leads to
\begin{eqnarray}
\nonumber \begin{array}{ll}
F_h^{ FEM}(t)=h\sum\limits_{j=1}^N\frac{\ddot v_{j+1} + 4 \ddot v_j+\ddot v_{j-1}}{6} j\frac{ v_{j+1}-v_{j-1}}{2}  +h\sum\limits_{j=1}^N\frac{\dot v_{j+1} + 4 \dot v_j+\dot v_{j-1}}{6} j\frac{ \dot v_{j+1}-\dot v_{j-1}}{2} &\\
 \quad +\frac{L}{6}\left(2\ddot v_{N+1}+\ddot v_{N}\right)(v_{N+1}-v_N)  +\frac{L}{6}\left(2\dot v_{N+1}+\dot v_{N}\right)(\dot v_{N+1}-\dot v_N)\\
 =-\frac{c^2 h}{2}\sum\limits_{j=0}^N \left|\frac{v_{j+1}-v_j}{h}\right|^2-\frac{h}{12}\sum\limits_{j=0}^N \left|\dot v_{j+1}+\dot v_j\right|^2 -\frac{h}{6}\sum\limits_{j=1}^N \dot v_{j+1}\dot v_j+\frac{c^2L}{2}\left|\frac{v_{N+1}-v_N}{h}\right|^2&\\
\quad -c^2 L\left|\frac{v_{N+1}-v_N}{h}\right|^2  -\frac{L\xi}{h} \dot v_{N+1}(v_{N+1}-v_N)+L\frac{(\dot v_{N+1}+\dot v_{N})^2}{12}+\frac{L}{6} \dot v_{N+1}\dot v_N&\\
\quad +L\frac{2\dot v_{N+1}^2-\dot v_{N+1}\dot v_N -\dot v_N^2}{6}&\\
=-E_h^{ FEM}(t) +\frac{h}{6}\sum\limits_{j=0}^N \left[\left|\dot v_j\right|^2-\dot v_{j+1}\dot v_j\right]-\frac{c^2 L}{2}\left|\frac{v_{N+1}-v_N}{h}\right|^2&\\
\quad  -\frac{L\xi}{h} \dot v_{N+1}(v_{N+1}-v_N) +L\frac{5\dot v_{N+1}^2+2\dot v_{N+1}\dot v_N -\dot v_N^2}{12}&\\
\le -E_h^{ FEM}(t) +\frac{h}{6}\sum\limits_{j=0}^N \left[\left|\dot v_j\right|^2-\dot v_{j+1}\dot v_j\right] +\frac{L}{2}(1+\frac{\xi^2}{c^2})|\dot v_{N+1}|^2.
 \end{array}
\end{eqnarray}
 \eqref{FEM-exp7} follows from this  together  with the following inequality
\begin{eqnarray}
\nonumber \begin{array}{ll}
\frac{h}{6}\sum\limits_{j=0}^N|\dot v_j|^2-\dot v_j \dot v_{j+1}&=\frac{h}{12}\sum\limits_{j=0}^N\frac{|\dot v_j - \dot v_{j+1}|^2}{2}-\frac{h}{12} |\dot v_{N+1}|^2\\
&\le  \frac{\kappa^{ FEM} h}{12c^2}\sum\limits_{j=0}^N|\dot v_j|^2\\
&\le \frac{\kappa^{ FEM}}{12c^2} E_h^{ FEM}(t).~~ \square
 \end{array}
\end{eqnarray}

\begin{theorem}\label{main2} Suppose that there exists a constant $0<\delta <\frac{2\xi c^2 }{L(c^2 +\xi^2)}$ such that for all initial conditions $\vec v^0, \vec v^1\in \mathbb{R}^{N+2},$ the energy $E_h^{ FEM}(t)$ corresponding to \eqref{FEM-exp}  satisfies $\forall t>0$
\begin{eqnarray}\label{thm-FEM}
E_h^{ FEM}(t)\le \frac{c+\delta L}{c-\delta L } e^{- \delta \left(1-\frac{L\delta}{c}\right)\left(1-\frac{\kappa^{ FEM}}{12c^2}\right)t} E_h^{ FEM}(0).
\end{eqnarray}
\end{theorem}
{\bf Proof:} Since $\frac{L_h^{ FEM}(t)}{dt}=\frac{dE^{ FEM}(t)}{dt}+\delta \frac{dF_h^{ FEM}(t)}{dt},$ by Lemma \ref{lem2}
\begin{eqnarray}
\nonumber \begin{array}{ll}
\frac{dL_h^{ FEM}(t)}{dt}&=-\left(\xi - \frac{L\delta}{2}\left(1+\frac{\xi^2}{c^2}\right)\right)|\dot v_{N+1}|^2 - \delta\left(1-\frac{\kappa^{ FEM} h^2}{12c^2}\right)E_h^{ FEM}(t)\\
&\le - \delta \left(1-\frac{L\delta}{c}\right)\left(1-\frac{\kappa^{ FEM}}{12c^2}\right)L^{FEM}_h(t).
 \end{array}
\end{eqnarray}
By the Gr\"onwall's inequality,
\begin{eqnarray}
\label{imp2}\begin{array}{ll}
L_h^{ FEM}(t)\le e^{-  \delta \left(1-\frac{L\delta}{c}\right)\left(1-\frac{\kappa^{ FEM}}{12c^2}\right)t}L_h^{ FEM}(0),
 \end{array}
\end{eqnarray}
and Lemma \ref{lem2}, \eqref{thm-FEM} is obtained. $\square$

\section{Maximal Decay Rate and Implementation of Direct Fourier Filtering}
For $N$ large enough (or small $h$) and $1\le k\le 2N+2$, Theorems \ref{thm-eigFD} and \ref{thm-eigFEM}  leads to
\begin{eqnarray}
\begin{array}{ll}\nonumber
\frac{h^2}{4c^2}\max \left\{\left|\lambda_k^{FD}(\xi,h)\right|^2\right\}&= \frac{\kappa^{FD}}{4c^2}=:\Gamma^{FD} + O(h),\\
\frac{h^2}{12c^2}\max \left\{\left|\lambda_k^{FEM}(\xi,h)\right|^2\right\}&=\frac{\kappa^{FEM}}{12c^2}=:\Gamma^{FEM} + O(h).
\end{array}
\end{eqnarray}
where $\Gamma$ is the Fourier filtering parameter. 
Therefore, for both model reductions $0<\Gamma^{FD}<1.$ 
Now, consider the space of filtered solutions for \eqref{eqr2}  and \eqref{eqr22} 
\begin{equation*}
	\mathcal{C}_h(\Gamma):=\left\{\vec{z}_h=\sum\limits_{0<\Gamma<1} a_k e^{i \sqrt{\lambda_k(\xi,h)}t}\vec{\phi}_k(h)\right\}.
\end{equation*}
By Theorems \ref{main1} and \ref{main2}, the exponential stability as  $h\to 0$ is immediate as some filtering $0<\Gamma<1$ is considered. If here is no filtering, i.e. $\Gamma\approx 1$, the exponential stability is at steak  since $1-\Gamma\to 0$ as $h\to 0$ e.g. see \eqref{imp1} and \eqref{imp2}. 

Notice that for each filtered solution in Theorems \ref{main1} and \ref{main2}, the decay rate $\sigma$ and $\delta$ are functions of $\Gamma$ and $\xi:$
\begin{eqnarray}\label{decayrate}
\begin{array}{ll}
\sigma (\Gamma, \xi)=\left\{\begin{array}{ll}
            \delta(1-\frac{L\delta}{c})(1-\Gamma^{FD}), & FD \\
        \delta \left(1-\frac{L\delta}{c}\right)\left(1-\Gamma^{FEM}\right), & FEM
    \end{array}\right.\\
   \delta (\xi)=\frac{c}{2L}\begin{array}{ll}
       \min \left( 1,  \frac{2\xi c }{c^2 +\xi^2}\right),&  FD~ {\text {and}} ~FEM. 
    \end{array}
    \end{array}
\end{eqnarray}
Note that $\delta$ reaches its maximum at $\xi=c:$
\begin{equation}\nonumber
\delta_{max} (\xi=c)=
     \frac{c}{2L}, \qquad  \text{FD ~ and ~ FEM} 
\end{equation}
at which $\sigma$ reaches its maximum
\begin{equation}\label{sigmamax}
\sigma_{max}\left(\delta_{max}\right)=
    \begin{cases}
        \frac{c}{4L}(1-\Gamma), \quad  \text{FD ~ and ~ FEM} 
    \end{cases}
\end{equation}
Our results perfectly mimic the maximal decay rate result \eqref{Lius} in Theorem \ref{thm1}.

\section{Numerical Experiments}

To show the strength of the Finite Difference \eqref{FD-exp} and Finite Element-based \eqref{FEM-exp} model reductions of  \eqref{estatic} with and without filtering, we consider high-enough number of nodes, e.g. $N=30$ with 60 complex eigenvalues of the system matrix in total.  Therefore, $h=1/31\approx 0.0322.$ For simplicity, consider $c=1$ and $L=1.$  For simulations, the following set of high-frequency initial conditions is considered
\begin{eqnarray}
&&\label{IC1}\left\{\begin{array}{ll}
 v(x_j,0)=v_j^0=10^{-3}\sum_{i=20}^{30} \sin(i \pi x_j)&\\
\dot  v(x_j,0)=v_j^1=10^{-3}\sum_{i=20}^{30} \sin(i \pi x_j).
\end{array}\right.
\end{eqnarray}

\begin{figure}[htb!]
    \centering
    \vspace{-0.02in}
              {{\includegraphics[width=5cm]{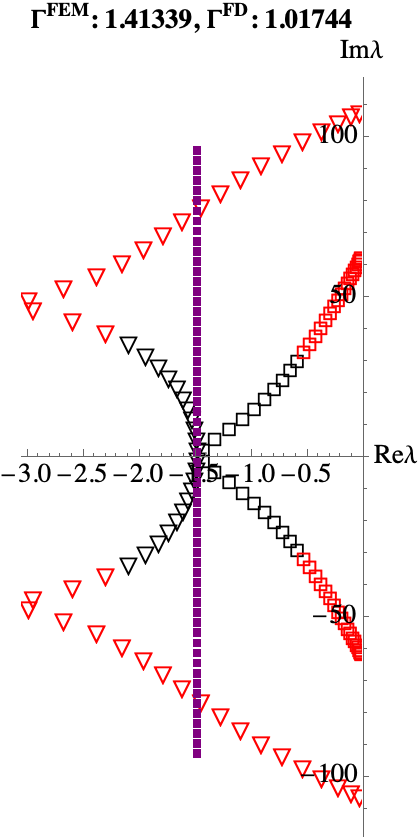} }}\\
             {{\includegraphics[width=8.8cm]{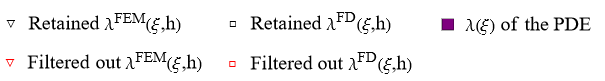} }}\\
        \vspace{0in}
    \caption{ For $c=1$ and $L=1$, 40 filtered out eigenvalues with $\xi=0.9<c=1$ for $2N=60$ eigenvalues. The filtered eigenvalues converge to the ones of the PDE.}%
    \label{spec}%
\end{figure}

The  eigenvalues $\lambda_k(\xi)$ of the PDE in \eqref{eigPDE} and  the approximated eigenvalues $\lambda_k^{FD}(\xi,h)$ and $\lambda_k^{FEM}(\xi,h)$ of the FEM and FD models, respectively, are simulated in Fig. \ref{spec}. For $\xi=0.9<c=1,$ the spectral plot shows   $40$ high-frequency eigenvalues out of  $60$ complex conjugate eigenvalues in total filtered out. The  corresponding   filtering parameters are $\Gamma^{FEM}=1.4133$ and $\Gamma^{FD}=1.017.$ Therefore, the maximal decay rates in \eqref{decayrate} are  calculated as $\sigma_{max}^{FEM}=0.2205$ and  $\sigma_{max}^{FD}=0.1864. $ This case can be compared to the maximal decay rate analysis in Theorem \ref{thm1} where the optimal decay rate is found to be $\sigma_{max}=0.25.$

\begin{figure}[htb!]
    \centering
    \vspace{-0.02in}
            {{\includegraphics[height=4.5cm]{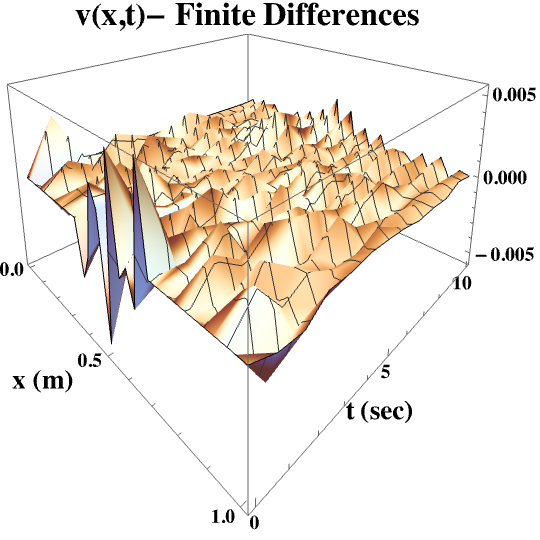} }}
                    {{\includegraphics[height=4.5cm]{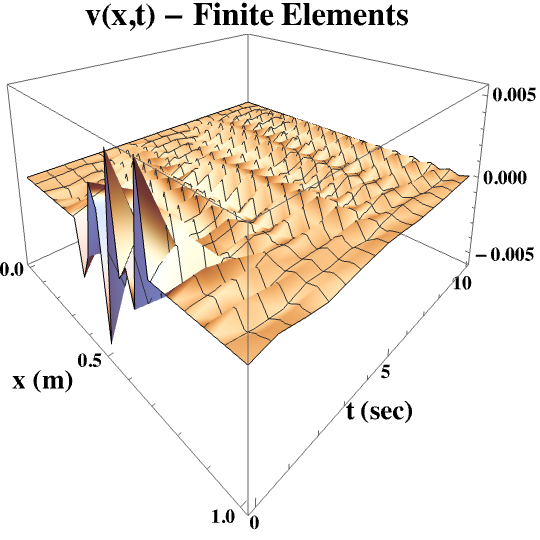} }}\\
                      {{\includegraphics[width=4.5cm]{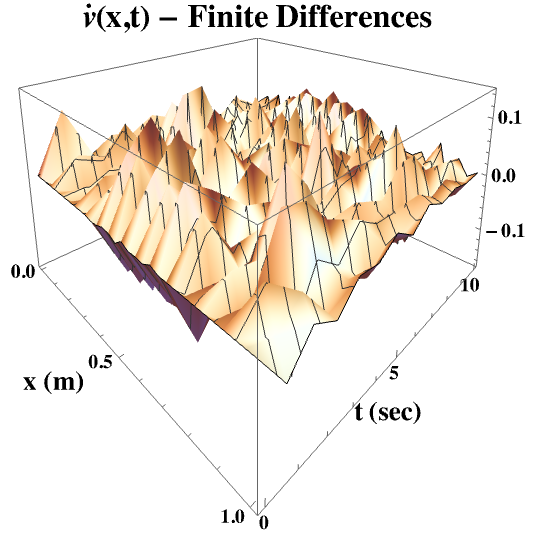} }}
            {{\includegraphics[width=4.5cm]{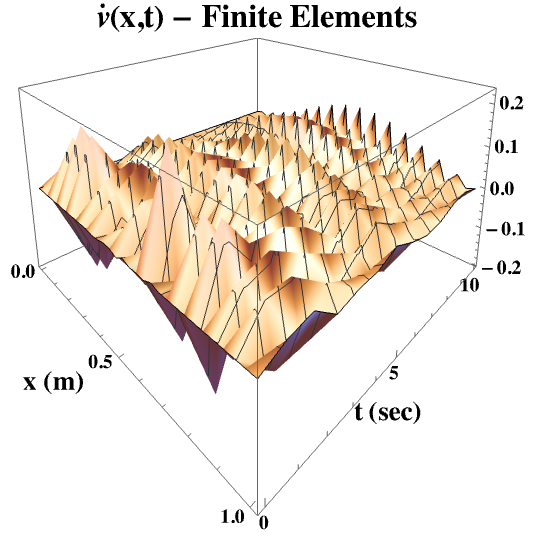} }}\\
                   {{\includegraphics[width=4.5cm]{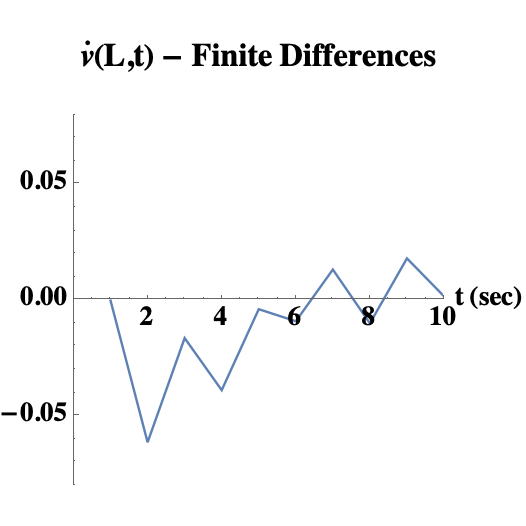} }}
            {{\includegraphics[width=4.5cm]{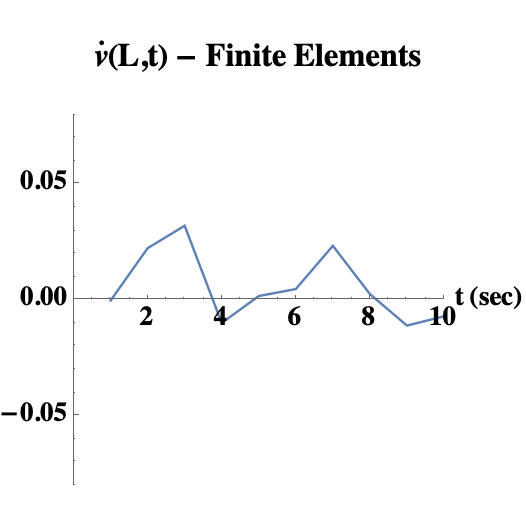} }}\\
        \vspace{0in}
    \caption{ If neither scheme is filtered,  the FEM model can be seen to stabilize the spurious high-frequency solutions much better than the FD model.  Yet, the solutions are still not stabilized to the rest position.}%
    \label{sim1}%
\end{figure}

\begin{figure}[htb!]
    \centering
    \vspace{-0.02in}
                   {{\includegraphics[height=4.5cm]{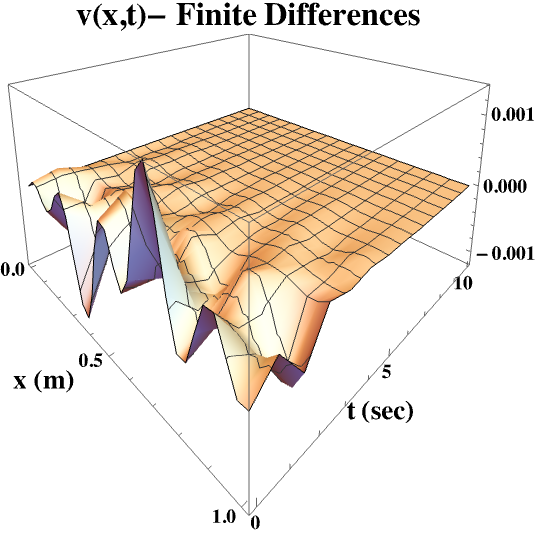} }}
                    {{\includegraphics[height=4.5cm]{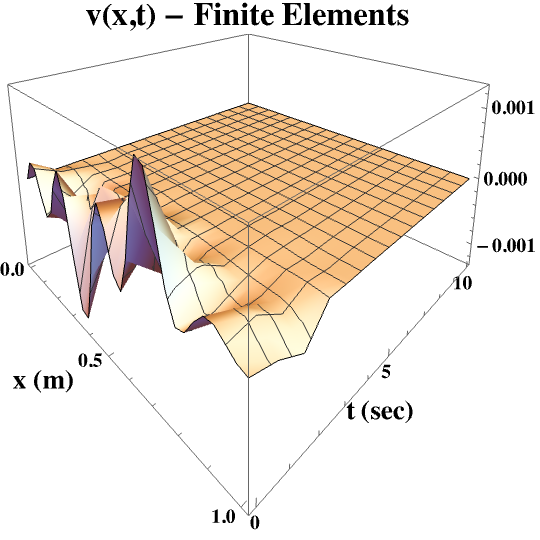} }}\\
                     {{\includegraphics[width=4.5cm]{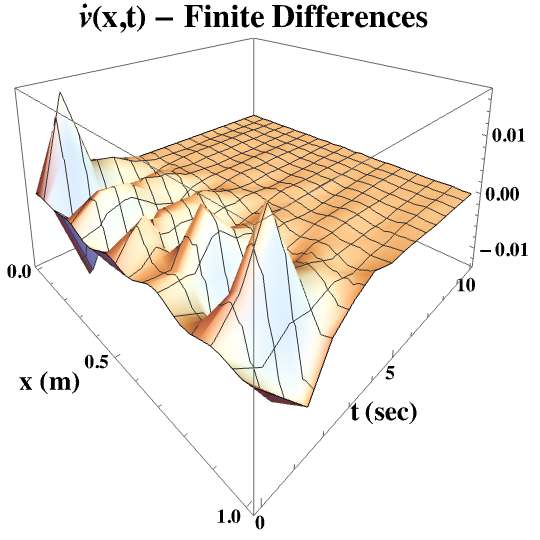} }}
            {{\includegraphics[width=4.5cm]{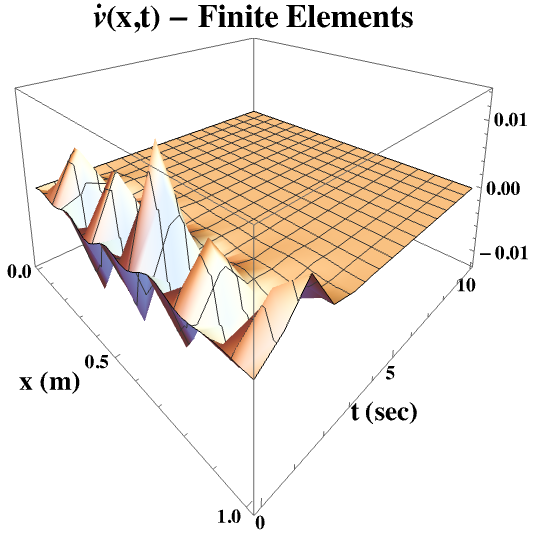} }}\\
                  {{\includegraphics[width=4.5cm]{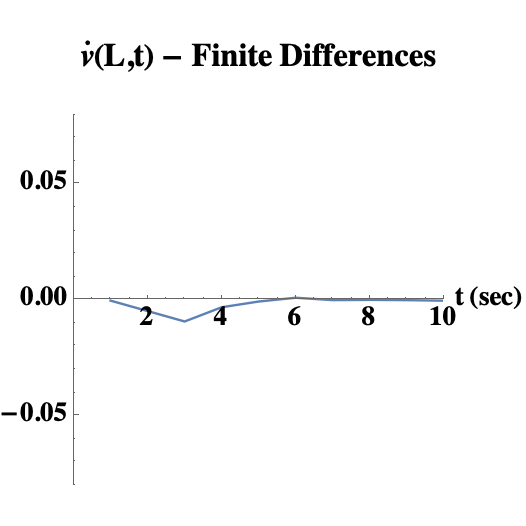} }}
            {{\includegraphics[width=4.5cm]{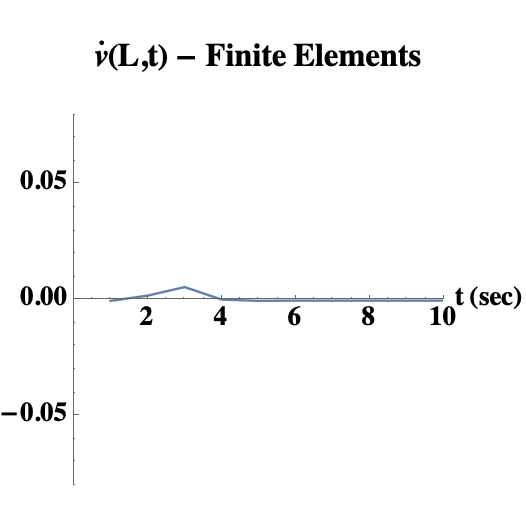} }}\\
        \vspace{0in}
    \caption{After filtering $40$ out of $60$ high-frequency solutions,  the FEM model stabilizes the high-frequency solutions much better than the FD model. }%
    \label{sim2}%
\end{figure}

Next, the simulations of $v(x,t), \dot v(x,t), E(t),$ and the tip velocity feedback $\dot v(L,t)$ with FD and FEM algorithms, $\xi=0.9<c=1,$ and for the unfiltered and filtered algorithms are shown in Figures \ref{sim1}-\ref{sim2}  for  the  set of initial conditions \eqref{IC1}, respectively.

Note that the simulations of various initial data (box, sawtooth, sinusoidal, pinch, square, triangle types) with and without filtering for both FD and FEM algorithms  are   real-time demonstrated in the recently-published Wolfram Demonstrations Project \cite{Wal}.

\section{Conclusions}
In conclusion, the Lyapunov approach laid out in this paper can be used to define an explicit decay rate in terms of the filtering parameter $\Gamma$ and the feedback gain $\xi.$ Our findings are in line with the conclusions in \cite{I-Z} that FEM provides a more accurate approximation and a better decay rate than FD to the \eqref{estatic}. For most PDE models based off of the wave \cite{M-O,Wilson} and beam equations \cite{LCSS,Guo2,Leon}, our approach can be easily adapted. Indeed, an immediate application under consideration is the obtension of the spectral estimates for coupled systems where there several branches of eigenvalues.

\end{document}